\documentclass{article}
\pdfoutput=1
\usepackage{amsmath}
\usepackage{amsfonts}
\usepackage{amssymb}
\usepackage{graphicx}%
\usepackage{tikz}
\usetikzlibrary{decorations.markings}
\usetikzlibrary{arrows}

\newtheorem{theorem}{Theorem}

\newtheorem{corollary}[theorem]{Corollary}

\newtheorem{definition}[theorem]{Definition}

\newtheorem{lemma}[theorem]{Lemma}

\newcommand{\GDP}{\mathcal{GDP}}
\newcommand{\join}{\textup{Join}}
\newenvironment{proof}[1][Proof]{\noindent\textbf{#1.} }{\ \rule{0.5em}{0.5em}}

\begin{document}

\title{A note on geometric duality in  matroid theory and knot theory}
\author{Lorenzo Traldi\\Lafayette College\\Easton Pennsylvania 18042}
\date{ }
\maketitle

\begin{abstract}
We observe that for planar graphs, the geometric duality relation generates both 2-isomorphism and abstract duality. This observation has the surprising consequence that for links, the equivalence relation defined by isomorphisms of checkerboard graphs is the same as the equivalence relation defined by 2-isomorphisms of checkerboard graphs.



\end{abstract}

\date{}

\section{Introduction}

About ninety years ago, Whitney published a series of papers including \cite{W1, W2, W3, W4}, in which he introduced the general theory of matroids, characterized graphs with isomorphic matroids, and related matroid duality to geometric duality of planar graphs. (For background and terminology we refer to Oxley's text \cite{Ox}, rather than Whitney's original papers.) Two of Whitney's particularly striking results are that a 3-connected loopless graph is characterized (up to isomorphism) by the associated cycle and cocycle matroids, and a 3-connected, loopless planar graph has a unique connected dual graph (abstract or geometric). In contrast, a planar graph $G$ that is not 3-connected may have different imbeddings on the 2-sphere $\mathbb S ^2$, with different geometric dual graphs. Moreover, if $G$ is not 2-connected it may have connected abstract dual graphs that are not geometric duals of any imbedding of $G$ on $\mathbb S ^2$.

In this note we observe that despite the latter fact, geometric duality suffices to define abstract duality of planar graphs. Moreover, geometric duality also suffices to define 2-isomorphism of planar graphs. 

\begin{theorem}
\label{main1}
Let $G$ and $G'$ be planar graphs. Then $G$ and $G'$ are abstract duals if and only if there exist an integer $k\geq 1$ and a sequence 
\[
G=G'_1,G'_2, \dots, G'_{2k}=G'
\]
such that $G'_i$ and $G'_{i+1}$ are geometric duals whenever $ 1 \leq i \leq 2k-1$.
\end{theorem}

\begin{theorem}
\label{main2}
Let $G$ and $G'$ be planar graphs. Then $G$ and $G'$ are 2-isomorphic if and only if there exist an integer $k\geq 1$ and a sequence 
\[
G=G'_1,G'_2, \dots, G'_{2k-1}=G'
\] 
such that $G'_i$ and $G'_{i+1}$ are geometric duals whenever $1 \leq i \leq 2k-2$.
\end{theorem}

Whitney's results imply that for 2-connected graphs, Theorem \ref{main1} can always be satisfied with $k=1$, and Theorem \ref{main2} can always be satisfied with $k=2$. 
In Sections 2 and 3 we deduce Theorems \ref{main1} and \ref{main2} from Whitney's theory. The arguments are not difficult; they require only some geometric manipulation of imbeddings of graphs with cut-vertices. We would not be surprised if Theorems \ref{main1} and \ref{main2} have appeared in the literature. If that is the case, we apologize for the failure of our repeated searches for an appropriate reference.

Theorems \ref{main1} and \ref{main2} may not seem to have much practical value, because it is easier to directly analyze the matroids associated with $G$ and $G'$ than it is to search through the set of all iterated geometric duals of all imbeddings of $G$ and $G'$ on $\mathbb S^2$. Nevertheless, Theorems \ref{main1} and \ref{main2} have some conceptual interest. In Section 4 we use them to prove the rather surprising result that for classical links, isomorphism of checkerboard graphs defines the same equivalence relation as 2-isomorphism of checkerboard graphs. 

\section{Geometric dual pairs}

In this note graphs may have loops and multiple edges. A \emph{planar} graph $G$ can be imbedded isomorphically in the plane $\mathbb R^2$, or equivalently the 2-sphere $\mathbb S^2$. To avoid topological pathologies we require the imbedded image of each edge to be a piecewise smooth curve of finite, positive length.

Here is a fairly standard definition.

\begin{definition}
\label{gdual}
Let $G$ be a planar graph. Then a \emph{geometric dual} of $G$ is a graph $H$ obtained (up to isomorphism) in the following way. First, imbed $G$ on $\mathbb S ^2$. Second, associate a vertex of $H$ with each complementary region of the imbedded image of $G$. Third, associate an edge of $H$ with each edge of $G$, so that the edge associated with $e \in E(G)$ is incident on each vertex of $H$ corresponding to a complementary region whose boundary includes the imbedded image of $e$. 
\end{definition}

It is easy to see that Definition \ref{gdual} is not symmetric: the primal graph $G$ can be any planar graph, but a geometric dual $H$ must be connected. (Geometric duality is symmetric for connected planar graphs, though.) When we want to indicate that Definition \ref{gdual} applies to two graphs, but we do not want to specify which graph is primal, we say simply that the two graphs are geometric duals. 

If $G$ and $H$ are geometric duals then the edge association of Definition \ref{gdual} defines a bijection $E(G) \leftrightarrow E(H)$. We call such a bijection a \emph{geometric duality map}.

Here are two other fairly standard definitions.

\begin{definition}
\label{adual}
Let $G$ and $H$ be graphs. An \emph{abstract duality map} between $G$ and $H$ is an edge bijection $E(G) \leftrightarrow E(H)$, under which edge sets of maximal forests of $G$ correspond to complements of edge sets of maximal forests of $H$.
\end{definition}

\begin{definition}
\label{2iso}
Let $G$ and $H$ be graphs. A \emph{2-isomorphism} between $G$ and $H$ is an edge bijection $E(G) \leftrightarrow E(H)$, under which edge sets of maximal forests of $G$ correspond to edge sets of maximal forests of $H$.
\end{definition}

Three fundamental properties of these definitions are: a geometric duality map is also an abstract duality map; the composition of two abstract duality maps is a 2-isomorphism; and the composition of an abstract duality map and a 2-isomorphism is an abstract duality map. We refer to Oxley \cite{Ox} for these properties, Whitney's results mentioned in the introduction, and other elements of the theory of graphs and matroids, including the partition of $E(G)$ as the disjoint union of the edge sets of the blocks of $G$. (A block of $G$ is either an isolated vertex or a connected subgraph with at least one edge, minimal with respect to the property of containing every circuit of which it contains an edge.)

\begin{definition}
\label{gdp}
Let $G$ and $H$ be connected planar graphs. Then the ordered pair $(G,H)$ is a \emph{geometric dual pair}, abbreviated \emph{gdp}, if $G$ and $H$ are geometric duals. The set of geometric dual pairs is denoted $\GDP$.
\end{definition}

\begin{definition}
\label{gdpsim} Let \emph{similarity}, denoted $\sim $, be the finest equivalence relation on $\mathcal{GDP}$
such that $(G,H) \sim (G',H')$ if $G$ is isomorphic to $G'$ or $H$ is isomorphic to $H'$.
\end{definition}

That is, $(G,H) \sim (G',H')$ if and only if there exist an integer $n \geq 1$ and a sequence 
\[
(G,H)=(G_1,H_1),(G_2,H_2), \dots, (G_n,H_n)=(G',H') \in \GDP
\]
such that whenever $1 \leq i \leq n-1$, $G_i \cong G_{i+1}$ or $H_i \cong H_{i+1}$. Notice that at each step of such a sequence we have a bijection $E(G_i) \to E(G_{i+1})$, given either by an isomorphism $G_i \cong G_{i+1}$ or by the composition of an isomorphism $H_i \cong H_{i+1}$ and geometric duality maps $E(G_i) \to E(H_i)$ and $E(H_{i+1}) \to E(G_{i+1})$. Either way, we have a 2-isomorphism $E(G_i) \to E(G_{i+1})$. The composition of these 2-isomorphisms is a 2-isomorphism between $G$ and $G'$, and the composition of this 2-isomorphism with geometric duality maps of $(G,H)$ and $(G',H')$ is a 2-isomorphism between $H$ and $H'$. We say these 2-isomorphisms are \emph{associated} with the similarity $(G,H) \sim (G',H')$.

We call a vertex $c \in V(G)$ a \emph{cut-vertex} of $G$ if $G-c$ has strictly more connected components than $G$ has, or $c$ is incident on a loop and another edge. That is, a cut-vertex is a vertex included in more than one block of $G$.

\begin{lemma}
\label{gdplem}
Suppose $(G,H) \in \GDP$ and $G$ has more than one cut-vertex. Then there is a gdp $(\widetilde G,\widetilde H) \in \GDP$ such that $(G,H) \sim (\widetilde G,\widetilde H)$, and $\widetilde G$ has strictly fewer cut-vertices than $G$ has.
\end{lemma}

\begin{proof}
If $G$ has cut-vertices $c \ne c'$, then $G$ can be imbedded on $\mathbb S ^2$ in the manner indicated schematically on the left in Figure \ref {figone}, with each of the subgraphs $G_1,G_2,G_3$ including at least one edge. ($G_1$ and $G_2$ both include $c$, and $G_2$ and $G_3$ both include $c'$.) We may choose the subgraphs $G_1,G_2,G_3$ so that $G_3$ contains every loop incident on $c'$, and every connected component of $G-c'$ other than the one containing $c$. Then $c'$ is not a cut-vertex of $G_2$.

Let $\widetilde G_3$ be the graph obtained from $G_3$ by replacing all adjacencies involving $c'$ with adjacencies involving $c$, and let $\widetilde G$ be the union of $G_1,G_2$ and $\widetilde G_3$. Then $\widetilde G$ may be imbedded on $\mathbb S ^2$ as indicated on the right in Figure \ref {figone}. 

\begin{figure}
\centering
\begin{tikzpicture} 
%left
\draw [thick, dotted] (-6,0) .. controls (-5.5,-0.5) and (-4.5,-0.5) .. (-4,0);
\draw [thick, dotted] (-6,0) .. controls (-5.5,0.5) and (-4.5,0.5) .. (-4,0);
\draw [thick, dotted] (-4,0) .. controls (-3.5,-0.5) and (-2.5,-0.5) .. (-2,0);
\draw [thick, dotted] (-4,0) .. controls (-3.5,0.5) and (-2.5,0.5) .. (-2,0);
\draw [thick, dotted] (-2,0) .. controls (-1.5,-0.5) and (-0.5,-0.5) .. (-0,0);
\draw [thick, dotted] (-2,0) .. controls (-1.5,0.5) and (-0.5,0.5) .. (-0,0);
%\draw [thick, dotted,  fill=white] (0,0) circle (1 cm);
\draw [thick, fill=black] (-4,0) circle (.1 cm);
\draw [thick, fill=black] (-2,0) circle (.1 cm);
\node at (-5,0) {$G_1$};
\node at (-3,0) {$G_2$};
\node at (-1,0) {$G_3$};
\node at (-4,0.5) {$c$};
\node at (-2,0.5) {$c'$};
%right
\draw [thick, dotted] (1.5,0) .. controls (2.0,-0.5) and (3.0,-0.5) .. (3.5,0);
\draw [thick, dotted] (1.5,0) .. controls (2.0,0.5) and (3.0,0.5) .. (3.5,0);
\draw [thick, dotted] (3.5,0) .. controls (4.0,-0.5) and (5.0,-0.5) .. (5.5,0);
\draw [thick, dotted] (3.5,0) .. controls (4.0,0.5) and (5.0,0.5) .. (5.5,0);
\draw [thick, dotted] (3.5,0) .. controls (3.0,-0.5) and (3.0,-1.5) .. (3.5,-2);
\draw [thick, dotted] (3.5,0) .. controls (4.0,-0.5) and (4.0,-1.5) .. (3.5,-2);
\draw [thick, fill=black] (3.5,0) circle (.1 cm);
\draw [thick, fill=black] (5.5,0) circle (.1 cm);
\node at (3.5,0.4) {$c$};
\node at (5.5,0.4) {$c'$};
\node at (2.5,0) {$G_1$};
\node at (4.5,0) {$G_2$};
\node at (3.5,-1) [rotate=-90] {$\widetilde G_3$};
\end{tikzpicture}
\caption{Imbeddings of $G$ and $\widetilde G$ in the proof of Lemma \ref{gdplem}.}
\label{figone}
\end{figure}
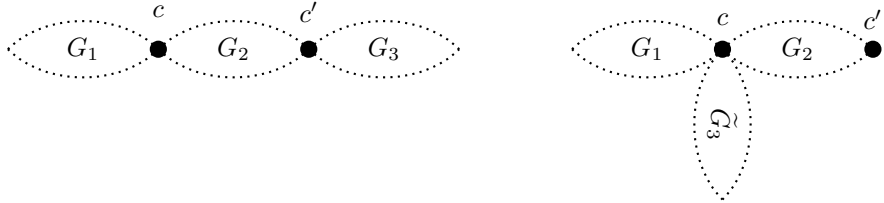

Let $\widetilde H$ be the geometric dual of the imbedding of $G$ indicated on the left in Figure \ref{figone}. Then $\widetilde H$ is obtained from the geometric duals of the imbeddings of $G_1,G_2$ and $G_3$ by identifying the three vertices corresponding to the region outside the dotted curves. The geometric dual of the imbedding of $\widetilde G$ indicated on the right in Figure \ref{figone} may be described in the same way, so it is isomorphic to $\widetilde H$. The lemma is satisfied because $(G,H) \sim (G,\widetilde H) \sim (\widetilde G,\widetilde H)$ and $\widetilde G$ has one fewer cut-vertex than $G$ has. 
\end{proof}

Here is the inductive generalization of Lemma \ref{gdplem}.

\begin{corollary}
\label{gdpcor}
If $G$ is a connected planar graph then there are gdps $(G,H) \sim (\widetilde G,\widetilde H)$ such that $\widetilde G$ does not have two different cut-vertices. 
\end{corollary}

\begin{definition}
\label{join}
Let $G_1,\dots,G_n$ be disjoint graphs, and suppose $v_i \in V(G_i)$ for $1 \leq i \leq n$. Then the \emph{one-point union} $\join(G_1, \dots, G_n; v_1, \dots, v_n)$ is the graph obtained from $\cup G_i$ by identifying $v_1, \dots, v_n$ to a single vertex.  
\end{definition}
\begin{lemma}
\label{gdplema}
Let $G_1$ and $G_2$ be disjoint, connected planar graphs. Let $G=\join(G_1, G_2; v_1, v_2)$ and $G'=\join(G_1, G_2; v'_1, v'_2)$, where $v_i,v'_i \in V(G_i)$. Then there are similar gdps $(G,H) \sim (G',H') \in \GDP$. 
\end{lemma}
\begin{proof}
Suppose first that $v_2=v'_2$. If $v_1$ and $v'_1$ are adjacent, then there are imbeddings of $G$ and $G'$ like the ones illustrated in Figure \ref{figtwo}, with the imbedded images of $G_2$ placed in the same complementary region of the imbedded image of $G_1$. These two imbeddings have the same geometric dual, $H$ say. Then $(G,H) \sim (G',H)$.

If $v_1$ and $v'_1$ are not adjacent, they are nonetheless connected by a path $v_1=w_1,w_2, \dots, w_k=v'_1$ in $G_1$. The observation of the preceding paragraph applies each time we walk from $w_i$ to $w_{i+1}$ on this path. As $\sim$ is an equivalence relation, we conclude that there are similar gdps $(G,H)$ and $(G',H')$.

\begin{figure}
\centering
\begin{tikzpicture} 
%left
\draw [thick] (-2,1.5) -- (-3,1.5) -- (-4,0.5) -- (-4,-0.5) -- (-3,-1.5) -- (-2,-1.5);
\draw [thick, dotted] (-4,0.5) .. controls (-3.5,1.0) and (-2.5,1.0) .. (-2,0.5);
\draw [thick, dotted] (-4,0.5) .. controls (-3.5,0.0) and (-2.5,0.0) .. (-2,0.5);
\draw [thick, fill=black] (-4,0.5) circle (.1 cm);
\draw [thick, fill=black] (-4,-0.5) circle (.1 cm);
\draw [thick, fill=black] (-3,1.5) circle (.1 cm);
\draw [thick, fill=black] (-3,-1.5) circle (.1 cm);
\draw [dotted, thick] (-2,1.5) .. controls (-1,0.5) and (-1,-0.5) .. (-2,-1.5);
\node at (-4,-1.5) {$G_1$};
\node at (-3,0.5) {$G_2$};
\node at (-4.8,0.5) {$v_1=v_2$};
\node at (-4.5,-0.5) {$v'_1$};
%right
\draw [thick] (6-2,1.5) -- (6-3,1.5) -- (6-4,0.5) -- (6-4,-0.5) -- (6-3,-1.5) -- (6-2,-1.5);
\draw [thick, dotted] (6-4,-0.5) .. controls (6-3.5,0) and (6-2.5,0) .. (6-2,-0.5);
\draw [thick, dotted] (6-4,-0.5) .. controls (6-3.5,-1) and (6-2.5,-1) .. (6-2,-0.5);
\draw [thick, fill=black] (6-4,0.5) circle (.1 cm);
\draw [thick, fill=black] (6-4,-0.5) circle (.1 cm);
\draw [thick, fill=black] (6-3,1.5) circle (.1 cm);
\draw [thick, fill=black] (6-3,-1.5) circle (.1 cm);
\draw [dotted, thick] (6-2,1.5) .. controls (6-1,0.5) and (6-1,-0.5) .. (6-2,-1.5);
\node at (6-4,-1.5) {$G_1$};
\node at (6-3,-0.5) {$G_2$};
\node at (6-4.9,-0.5) {$v'_1=v_2$};
\node at (6-4.5,0.5) {$v_1$};
\end{tikzpicture}
\caption{Imbeddings of $G$ and $G'$ in the proof of Lemma \ref{gdplema}.}
\label{figtwo}
\end{figure}
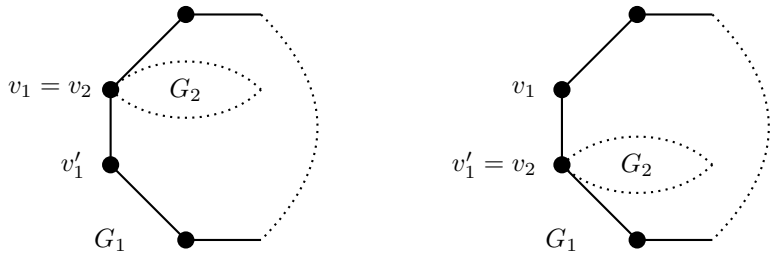

If $v_1=v'_1$, then the same argument applies, with $G_1$ and $G_2$ interchanged. The general case follows.
\end{proof}

Here is the inductive generalization of Lemma \ref{gdplema}.
\begin{corollary}
\label{gdpcora}
Let $G_1, \dots, G_n$ be disjoint, connected planar graphs. Suppose $G=\join(G_1, \dots, G_n; v_1, \dots, v_n)$ and $G'=\join(G_1, \dots, G_n; v'_1, \dots, v'_n)$. Then there are similar gdps $(G,H) \sim (G',H') \in \GDP$.  
\end{corollary}

Combining Corollaries \ref{gdpcor} and \ref{gdpcora}, we obtain the following.

\begin{corollary}
\label{gdpcorb}
Let $G$ be a connected planar graph, and let $B_1, \dots, B_n$ be the blocks of $G$, considered as disjoint graphs. For each $i \in \{1, \dots, n \}$, let $v_i$ be a vertex of $B_i$. Then there is a similarity $(G,H) \sim (\join(B_1, \dots, B_n; v_1, \dots, v_n), \widetilde H)$, with an associated 2-isomorphism that is the identity map on $E(B_i)$, for each $i$.
\end{corollary}
\begin{proof}
The similarities of Lemmas \ref{gdplem} and \ref{gdplema} involve detaching and re-attaching subgraphs at different cut-vertices, without disturbing the internal structure of any block.
\end{proof}
\begin{theorem} 
\label{gdpthm}
Every 2-isomorphism between connected planar graphs $G$ and $G'$ is associated with a similarity $(G,H) \sim (G',H')$.
\end{theorem}

\begin{proof}
Let $B_1, \dots, B_n$ be the blocks of $G$, considered as disjoint graphs. Let $v_i \in V(B_i)$ for each $i$, and let $\widetilde G = \join(B_1, \dots, B_n; v_1, \dots, v_n)$. Corollary \ref{gdpcorb} gives us a similarity $(G,H) \sim (\widetilde G, \widetilde H)$ with an associated edge bijection that is the identity map on $E(B_i)$, for each $i$.

$\widetilde G$ has an imbedding on $\mathbb S ^2$ of the type depicted in Figure \ref{figthree}. For $1 \leq i \leq n$, let $H_i$ be the geometric dual of the imbedding of $B_i$ that appears in Figure \ref{figthree}. Then the geometric dual of Figure \ref{figthree} is $\widetilde H=\join(H_1, \dots, H_n; w_1, \dots, w_n)$, where $w_i$ always corresponds to the complementary region outside the dotted lines in Figure \ref{figthree}.

Suppose $f:E(G) \to E(G')$ is a 2-isomorphism. Then the blocks of $G'$ are the connected subgraphs $B'_1, \dots, B'_n$ with edge sets $f(E(B_1)), \dots, f(E(B_n))$. For each index $i$, the restriction of $f$ to $B_i$ is a 2-isomorphism with $B'_i$, and hence defines an abstract duality between $H_i$ and $B'_i$. As $B'_i$ is a block, $H_i$ can be re-imbedded in $\mathbb S^2$ so as to be geometrically dual to $B'_i$, with $f$ the geometric duality bijection. (See \cite[Prop.\ 5.2.3]{Ox} or \cite[Thm.\ 30]{W2}; the map $f$ is implicit in those discussions.) Combining disjoint copies of these imbeddings, we obtain an imbedding of a one-point union $\widetilde G' = \join(B'_1, \dots, B'_n; v'_1, \dots, v'_n)$ that resembles Figure \ref{figthree}, with each $B_i$ replaced by $B'_i$. Therefore $(\widetilde G,\widetilde H) \sim (\widetilde G',\widetilde H)$. 

\begin{figure}
\centering
\begin{tikzpicture} 
\draw [thick, dotted] (1.5,0) .. controls (2.0,-0.4) and (3.0,-0.4) .. (3.5,0);
\draw [thick, dotted] (1.5,0) .. controls (2.0,0.4) and (3.0,0.4) .. (3.5,0);

\draw [thick, dotted] (3.5,0) .. controls (2.5,0.5) and (2.3,1.5) .. (2.3,1.8);
\draw [thick, dotted] (3.5,0) .. controls (3.4,0.5) and (3.3,1.5) .. (2.3,1.8);

\draw [thick, dotted] (3.5,0) .. controls (3.6,0.5) and (3.7,1.5) .. (4.7,1.8);
\draw [thick, dotted] (3.5,0) .. controls (4.5,0.5) and (4.7,1.5) .. (4.7,1.8);

\draw [thick, dotted] (3.5,0) .. controls (2.5,-0.5) and (2.3,-1.5) .. (2.3,-1.8);
\draw [thick, dotted] (3.5,0) .. controls (3.4,-0.5) and (3.3,-1.5) .. (2.3,-1.8);
\draw [thick, fill=black] (3.5,0) circle (.1 cm);
\draw [thick, dashed, domain=-90:30]  plot ({3.5+(1.2)*cos(\x)}, {(1.2)*sin(\x)});
\node at (2.5,0) {$B_1$};
\node at (2.85,-1) {$B_n$};
\node at (2.9,1) {$B_2$};
\node at (4.1,1) {$B_3$};
\end{tikzpicture}
\caption{An imbedding of $\widetilde G$ in the proof of Theorem \ref{gdpthm}.}
\label{figthree}
\end{figure}
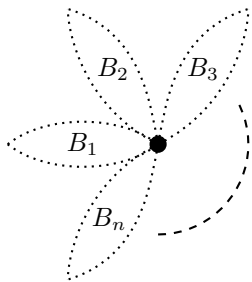

According to Corollary \ref{gdpcorb}, there is a similarity $(\widetilde G',\widetilde H) \sim (G', H')$ for some geometric dual $H'$ of $G'$, with an associated 2-isomorphism that is the identity map on each $E(B'_i)$. Composing the similarities $(G,H) \sim (\widetilde G, \widetilde H) \sim (\widetilde G', \widetilde H) \sim (G',H')$ gives us a similarity $(G,H) \sim (G',H')$, with an associated edge bijection that agrees with $f$ on every block of $G$, and hence must be $f$.
\end{proof}

\section{Theorems \ref{main1} and \ref{main2}}

The proofs of Theorems \ref{main1} and \ref{main2} are almost completely covered by Whitney's theory and the argument of the preceding section. We provide the few remaining details here. 

We begin with Theorem \ref{main2}. Replacing each of $G,G'$ by a geometric dual of a geometric dual, we may assume without loss of generality that $G$ and $G'$ are connected. The ``if'' direction of Theorem \ref{main2} then follows from the fact that the composition of two geometric duality maps is a 2-isomorphism.

For the ``only if'' direction of Theorem \ref{main2}, suppose $G$ and $G'$ are 2-isomorphic planar graphs. According to Theorem \ref{gdpthm}, there are gdps $(G,H) \sim (G',H') \in \GDP$. That is, there is a sequence 
\begin{equation}
\label{gdpseq}
(G,H)=(G_1,H_1), (G_2,H_2), \dots, (G_n,H_n)=(G',H') \in \GDP
\end{equation}
such that for each $i \in \{1, \dots, n-1\}$, either $G_i \cong G_{i+1}$ or $H_i \cong H_{i+1}$. Choose such a sequence with $n$ as small as possible. If $n=1$ then $G \cong G'$ and $H \cong H'$, so Theorem \ref{main2} is satisfied by the sequence $G,H,G'$.

Suppose $n>1$. Notice that we may remove the gdp $(G_{i+1},H_{i+1})$ from (\ref{gdpseq}) if $G_i \cong G_{i+1} \cong G_{i+2}$ or $H_i \cong H_{i+1} \cong H_{i+2}$. Similarly, if $G_1 \cong G_2$ we may remove $(G_1,H_1)$, and if $G_{n-1} \cong G_n$ we may remove $(G_n,H_n)$. As $n$ is minimal, none of these removals is possible. Therefore $G_i \cong G_{i+1}$ when $i$ is even, $H_i \cong H_{i+1}$ when $i$ is odd, and $n$ is even. Hence Theorem \ref{main2} is satisfied by the sequence 
\[
G=G_1,H_2,G_2,H_4,G_4, \dots, H_n, G_n = G'.
\]

For Theorem \ref{main1}, let $G$ and $G'$ be planar graphs, and let $H$ be a geometric dual of $G$. Then $G$ and $G'$ are abstract duals if and only if $H$ and $G'$ are 2-isomorphic, so we may deduce Theorem \ref{main1} for $G$ and $G'$ by applying Theorem \ref{main2} to $H$ and $G'$.

Notice that Theorem \ref{gdpthm} allows us to refine Theorems \ref{main1} and \ref{main2} to results about maps rather than graphs: every 2-isomorphism between planar graphs is a composition of an even number of geometric duality maps, and every abstract duality map between planar graphs is a composition of an odd number of geometric duality maps. These results will be useful in the next section.

\section{An application to knot theory}

We begin this section with a brief review of some basic ideas of classical knot theory. Thorough treatments are presented in many textbooks, like \cite{Li}.

A (classical) \emph{link} is a collection of finitely many, pairwise disjoint, piecewise smooth simple closed curves in three-space $\mathbb R^3$ (or equivalently, the three-sphere $\mathbb S^3$). The curves that constitute a link are its \emph{components}. Two links are \emph{equivalent} if one can be transformed into the other by continuously moving it around in space. A \emph{diagram} of a link is obtained from a regular projection -- that is, a projection in a plane with only finitely many singularities, all of them transverse double points called \emph{crossings} -- by removing two short arcs at each crossing, to indicate which of the segments is farther from the observer. It is a famous fact that diagrams of equivalent links are connected to each other by three types of diagram changes, the Reidemeister moves.

The complementary regions of a link diagram can be 2-colored, checkerboard fashion; we call the colors ``shaded'' and ``unshaded.'' Such a 2-coloring defines a pair of graphs associated with the link diagram; one graph has vertices corresponding to the shaded regions, and the other has vertices corresponding to the unshaded regions. Both graphs have edges corresponding to crossings. 

\begin{figure} 
\centering
\begin{tikzpicture}
%leftlink
\draw [thick] [white, fill=lightgray!50] (-.6,1.15) -- (-.6,1.45) -- (-.9,1.45) -- (-.9,1.6) -- (2.4,1.6) -- (2.4,1) -- (2,1) -- (2,1.3) -- (.2,1.3) -- (.2,1) -- (-.6,1) -- (-.6,1.15);
\draw [thick] [white, fill=lightgray!50] (-.6,.8) -- (-.6,1) -- (-.9,1) -- (-.9,1.45) -- (-1.6,1.45) -- (-1.6,.8) -- (-1.45,.8) -- (-1.45,1.2) -- (-1.15,1.2) -- (-1.15,.8);
\draw [thick] [white, fill=lightgray!50] (.2,1) -- (2,1) -- (2,.8) -- (1.7,.5) -- (1.4,.8) -- (1.1,.5) -- (.8,.8) -- (.5,.5) -- (.2,.8);
\draw [thick] [white, fill=lightgray!50] (.2,0) -- (2,0) -- (2,.2) -- (1.7,.5) -- (1.4,.2) -- (1.1,.5) -- (.8,.2) -- (.5,.5) -- (.2,.2);
\draw [thick] [white, fill=lightgray!50] (-1.45,.8) -- (-1.15,.8) -- (-1.15,.5) -- (-1.45,.5);
%\draw [thick, fill=lightgray!50] (-1,.8) circle (.3 cm);
%\draw [thick, fill=white] (-1,.8) circle (.15 cm);
\draw [thick] (0.2,0) -- (2,0);
\draw [thick] (0.2,0) -- (0.2,0.2);
\draw [thick] (2,.2) -- (2,0);
\draw [thick] (.2,.2) -- (.4,.4);
\draw [thick] (.6,.6) -- (.8,.8);
\draw [thick] (.8,.8) -- (1.4,.2);
\draw [thick] (.2,.8) -- (.8,0.2);
\draw [thick] (.8,.2) -- (1,0.4);
\draw [thick] (1.4,.8) -- (1.2,0.6);
\draw [thick] (1.4,.8) -- (2,0.2);
\draw [thick] (.8,.2) -- (1,0.4);
\draw [thick] (1.4,.2) -- (1.6,0.4);
\draw [thick] (1.8,.6) -- (2,0.8);
\draw [thick] (2,.9) -- (2,0.8);
\draw [thick] (0.2,1.3) -- (0.2,0.8);
\draw [thick] (0.2,1.3) -- (2,1.3);
\draw [thick] (2,1.1) -- (2,1.3);
\draw [thick] (-.9,1) -- (.1,1);
\draw [thick] (0.3,1) -- (2.4,1);
\draw [thick] (-.9,1) -- (-.9,1.6);
\draw [thick] (2.4,1.6) -- (-.9,1.6);
\draw [thick] (2.4,1.6) -- (2.4,1);
%\draw [thick] (-1,1.45) -- (-1.6,1.45) -- (-1.6,1.15) -- (-1,1.15);
%\draw [thick] (-.6,1.45) -- (-.8,1.45) -- (-.8,1.15) -- (-.6,1.15);
\draw [thick] (-1,1.45) -- (-1.6,1.45) -- (-1.6,.8) -- (-.6,.8) -- (-.6,.9);
\draw [thick] (-.8,1.45) -- (-.6,1.45) -- (-.6,1.1);
\draw [thick] (-1.15,.9) -- (-1.15,1.2) -- (-1.45,1.2) -- (-1.45,.9);
\draw [thick] (-1.15,.7) -- (-1.15,.5) -- (-1.45,.5) -- (-1.45,.7);
%leftgraph
\draw [thick, fill=black] (1+3.7,1.5) circle (.1 cm);
\draw [thick, fill=black] (1+3.7,.8) circle (.1 cm);
\draw [thick] (1+4.5,1.5) .. controls (1+4.15,1.7)  .. (1+3.7,1.5);
\draw [thick] (1+4.5,1.5) .. controls (1+4.15,1.3)  .. (1+3.7,1.5);
\draw [thick] (1+3.7,1.5) .. controls (1+3.5,1.15)  .. (1+3.7,.8);
\draw [thick] (1+3.7,1.5) .. controls (1+3.9,1.15)  .. (1+3.7,.8);
\draw [thick] (1+4.5,0.1) .. controls (1+4.2,.45)  .. (1+4.5,.8);
\draw [thick] (1+4.5,0.1) .. controls (1+4.8,.45)  .. (1+4.5,.8);
\draw [thick] (1+4.5,0.1) -- (1+4.5,0.8);
\draw [thick] (1+4.5,1.5) .. controls (1+4.25,1.15)  .. (1+4.5,.8);
\draw [thick] (1+4.5,1.5) .. controls (1+4.75,1.15)  .. (1+4.5,.8);
\draw [thick, fill=black] (1+4.5,0.1) circle (.1 cm);
\draw [thick, fill=black] (1+4.5,0.8) circle (.1 cm);
\draw [thick, fill=black] (1+4.5,1.5) circle (.1 cm);
%rightgraph
\draw [thick, fill=black] (2+5.7,1.1) circle (.1 cm);
\draw [thick, fill=black] (2+6.2,1.3) circle (.1 cm);
\draw [thick] (2+7.2,-0.1) .. controls (2+6.2,.5)  .. (2+7.2,1.1);
\draw [thick] (2+7.2,-0.1) .. controls (2+8.2,.5)  .. (2+7.2,1.1);
\draw [thick] (2+7.2,-0.1) .. controls (2+5.6,0)  .. (2+5.7,1.1);
\draw [thick] (2+7.2,-0.1) .. controls  (2+5.9,0.1)  .. (2+5.7,1.1);
\draw [thick] (2+7.2,-0.1) .. controls (2+6,0.2)  .. (2+6.2,1.3);
\draw [thick] (2+7.2,-0.1) .. controls  (2+6.3,0.3)  .. (2+6.2,1.3);
\draw [thick] (2+6.9,0.5) -- (2+7.5,0.5);
\draw [thick] (2+6.9,0.5) -- (2+7.2,-.1);
\draw [thick] (2+7.5,0.5) -- (2+7.2,-.1);
\draw [thick, fill=black] (2+7.2,-.1) circle (.1 cm);
\draw [thick, fill=black] (2+6.9,0.5) circle (.1 cm);
\draw [thick, fill=black] (2+7.5,0.5) circle (.1 cm);
\draw [thick, fill=black] (2+7.2,1.1) circle (.1 cm);
\end{tikzpicture}
\caption{An equivalent diagram, with connected checkerboard graphs.}
\label{figfive}
\end{figure}
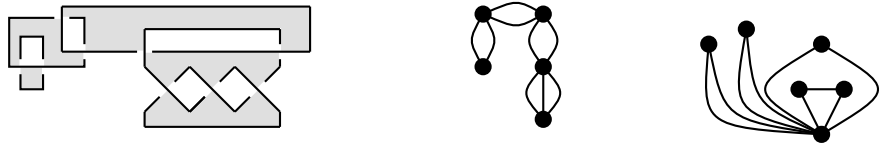

These checkerboard graphs are geometric duals if one is connected, but they can both be disconnected, as in Figure \ref{figfour}. Figure \ref{figfive} illustrates a simple trick used to restrict attention to link diagrams for which both checkerboard graphs are connected. Assuming both are connected is convenient because it guarantees that valuable knot-theoretic information is incorporated in both checkerboard graphs. For example, if we remove one of the unknotted link components represented by concentric squares on the left in Figure \ref{figfour}, we obtain a new diagram $D'$, which shares a disconnected checkerboard graph with the diagram pictured in Figure \ref{figfour}. It is not hard to show that the link represented by $D'$ does not have a dual pair of connected checkerboard graphs either of which is shared with the link from Figures \ref{figfour} and \ref{figfive}.

%For example, notice that each of the unknotted link components represented by concentric squares on the left in Figure \ref{figfour} is only represented in one checkerboard graph; in contrast, both of them are represented in both checkerboard graphs in Figure \ref{figfive}.

\begin{figure} 
\centering
\begin{tikzpicture}
%leftlink
\draw [thick] [white, fill=lightgray!50] (-.2,1.6) -- (2.4,1.6) -- (2.4,1) -- (2,1) -- (2,1.3) -- (.2,1.3) -- (.2,1) -- (-.2,1);
\draw [thick] [white, fill=lightgray!50] (.2,1) -- (2,1) -- (2,.8) -- (1.7,.5) -- (1.4,.8) -- (1.1,.5) -- (.8,.8) -- (.5,.5) -- (.2,.8);
\draw [thick] [white, fill=lightgray!50] (.2,0) -- (2,0) -- (2,.2) -- (1.7,.5) -- (1.4,.2) -- (1.1,.5) -- (.8,.2) -- (.5,.5) -- (.2,.2);
\draw [thick] [white, fill=lightgray!50] (-.95,.8) -- (-.95,1.5) -- (-1.65,1.5) -- (-1.65,.8) -- (-.95,.8);
\draw [thick] [white, fill=white] (-1.15,1) -- (-1.15,1.3) -- (-1.45,1.3) -- (-1.45,1) -- (-1.15,1);
%\draw [thick, fill=lightgray!50] (-1,.8) circle (.3 cm);
%\draw [thick, fill=white] (-1,.8) circle (.15 cm);
\draw [thick] (0.2,0) -- (2,0);
\draw [thick] (0.2,0) -- (0.2,0.2);
\draw [thick] (2,.2) -- (2,0);
\draw [thick] (.2,.2) -- (.4,.4);
\draw [thick] (.6,.6) -- (.8,.8);
\draw [thick] (.8,.8) -- (1.4,.2);
\draw [thick] (.2,.8) -- (.8,0.2);
\draw [thick] (.8,.2) -- (1,0.4);
\draw [thick] (1.4,.8) -- (1.2,0.6);
\draw [thick] (1.4,.8) -- (2,0.2);
\draw [thick] (.8,.2) -- (1,0.4);
\draw [thick] (1.4,.2) -- (1.6,0.4);
\draw [thick] (1.8,.6) -- (2,0.8);
\draw [thick] (2,.9) -- (2,0.8);
\draw [thick] (0.2,1.3) -- (0.2,0.8);
\draw [thick] (0.2,1.3) -- (2,1.3);
\draw [thick] (2,1.1) -- (2,1.3);
\draw [thick] (-.2,1) -- (.1,1);
\draw [thick] (0.3,1) -- (2.4,1);
\draw [thick] (-0.2,1) -- (-.2,1.6);
\draw [thick] (2.4,1.6) -- (-.2,1.6);
\draw [thick] (2.4,1.6) -- (2.4,1);
\draw [thick] (-.95,.8) -- (-.95,1.5) -- (-1.65,1.5) -- (-1.65,.8) -- (-.95,.8);
\draw [thick] (-1.15,1) -- (-1.15,1.3) -- (-1.45,1.3) -- (-1.45,1) -- (-1.15,1);
%leftgraph
\draw [thick, fill=black] (1+3.8,1.15) circle (.1 cm);
\draw [thick] (1+4.5,0.1) .. controls (1+4.2,.45)  .. (1+4.5,.8);
\draw [thick] (1+4.5,0.1) .. controls (1+4.8,.45)  .. (1+4.5,.8);
\draw [thick] (1+4.5,0.1) -- (1+4.5,0.8);
\draw [thick] (1+4.5,1.5) .. controls (1+4.2,1.15)  .. (1+4.5,.8);
\draw [thick] (1+4.5,1.5) .. controls (1+4.8,1.15)  .. (1+4.5,.8);
\draw [thick, fill=black] (1+4.5,0.1) circle (.1 cm);
\draw [thick, fill=black] (1+4.5,0.8) circle (.1 cm);
\draw [thick, fill=black] (1+4.5,1.5) circle (.1 cm);
%rightgraph
%\draw [thick, fill=black] (6,.8) circle (.1 cm);
\draw [thick, fill=black] (2+6.2,1.15) circle (.1 cm);
\draw [thick] (2+7.2,-0.1) .. controls (2+6.2,.5)  .. (2+7.2,1.1);
\draw [thick] (2+7.2,-0.1) .. controls (2+8.2,.5)  .. (2+7.2,1.1);
\draw [thick] (2+6.9,0.5) -- (2+7.5,0.5);
\draw [thick] (2+6.9,0.5) -- (2+7.2,-.1);
\draw [thick] (2+7.5,0.5) -- (2+7.2,-.1);
\draw [thick, fill=black] (2+7.2,-.1) circle (.1 cm);
\draw [thick, fill=black] (2+6.9,0.5) circle (.1 cm);
\draw [thick, fill=black] (2+7.5,0.5) circle (.1 cm);
\draw [thick, fill=black] (2+7.2,1.1) circle (.1 cm);
\end{tikzpicture}
\caption{A link diagram, with its checkerboard graphs.}
\label{figfour}
\end{figure}
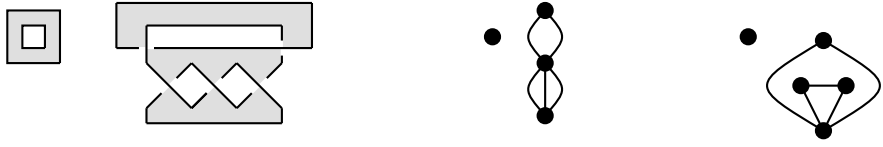

\begin{figure}
\centering
\begin{tikzpicture} 
\draw [thick] [white, fill=lightgray!50] (-2,.5) -- (-1.5,0) -- (-1,.5);
\draw [thick] [white, fill=lightgray!50] (-2,-.5) -- (-1.5,0) -- (-1,-.5);
\draw [thick] [white, fill=lightgray!50] (2,.5) -- (1.5,0) -- (1,.5);
\draw [thick] [white, fill=lightgray!50] (2,-.5) -- (1.5,0) -- (1,-.5);
\draw [thick] (-2,.5) -- (-1,-.5);
\draw [thick] (-2,-.5) -- (-1.6,-.1);
\draw [thick] (-1.4,.1) -- (-1,.5);
\draw [thick] (2,-.5) -- (1.6,-.1);
\draw [thick] (1,-.5) -- (2,.5);
\draw [thick] (1.4,.1) -- (1,.5);
\node at (-2.5,0) {$-1$};
\node at (2.5,0) {$1$};
\end{tikzpicture}
\caption{Crossing signs used as edge signs in the unshaded checkerboard graph.}
\label{figsix}
\end{figure}
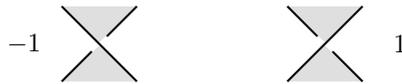

The edges of an unshaded checkerboard graph are given $\pm$ signs using the convention of Figure \ref{figsix}, and the edges of the shaded checkerboard graph are given the opposite signs. Using opposite signs has the effect of guaranteeing that if the 2-coloring of the complementary regions is reversed, the two signed checkerboard graphs are not affected; they are simply interchanged.

Here is a natural definition.

\begin{definition}
\label{isosign}
An isomorphism (resp., 2-isomorphism) of signed graphs $G$ and $H$ is an isomorphism (resp., 2-isomorphism) $f$ between the unsigned versions of $G$ and $H$, such that for every $e \in E(G)$, $e$ and $f(e)$ have the same sign. 
\end{definition}

For general graphs, isomorphism is a much stricter requirement than 2-isomorphism. Despite this difference, it turns out that for link diagrams, isomorphism and 2-isomorphism of signed checkerboard graphs define the same equivalence relation.

\begin{theorem}
\label{main3} Let $\mathcal{D}$ be the set of link diagrams with both checkerboard graphs connected.
Let $\sim_1$ be the finest equivalence relation on $\mathcal{D}$ with $D \sim_1 D'$ whenever one of the signed checkerboard graphs of $D$ is isomorphic to one of the signed checkerboard graphs of $D'$. Let $\sim_2$ be the finest equivalence relation on $\mathcal{D}$ with $D \sim_2 D'$ whenever one of the signed checkerboard graphs of $D$ is 2-isomorphic to one of the signed checkerboard graphs of $D'$. Then $\sim_1$ and $\sim_2$ are the same equivalence relation on $\mathcal{D}$.
\end{theorem}

\begin{proof} If $D \sim_1 D'$ then it is clear that $D \sim_2 D'$ too.

For the converse, note first that if $(G,H) \in \GDP$ and corresponding edges in $G$ and $H$ are given opposite signs, then it is possible to construct a link diagram for which these signed versions of $G$ and $H$ are the signed checkerboard graphs. Simply imbed $G$ and $H$ as geometric duals in the plane, form the associated medial graph, and choose the crossings to match the edge signs.

Now, suppose $D$ and $D'$ are link diagrams with connected checkerboard graphs, and there is a 2-isomorphism $f:E(G) \to E(G')$ where $G$ is a signed checkerboard graph of $D$ and $G'$ is a signed checkerboard graph of $D'$. According to the refined version of Theorem \ref{main2} mentioned at the end of Section 3, $f$ is a composition of an even number of geometric duality maps. That is, there exist a positive integer $k\geq 1$, connected planar graphs $G=G_1, \dots, G_{2k-1}=G'$, and  geometric duality maps $f_i:E(G_i) \to E(G_{i+1})$ with $f= f_{2k-2} \circ f_{2k-3} \circ \cdots \circ f_1$. If $k=1$ then $G=G'$, so $D \sim_1 D'$.

Suppose $k>1$. Define edge signs on $G_2, \dots, G_{2k-1}$ so that if $1 \leq i \leq 2k-2$ and $e \in E(G_i)$, the sign of $f_i(e)$ is the opposite of the sign of $e$. Then the resulting edge signs on $G_{2k-1}$ are the same as the original edge signs of $G'$. For $1 \leq i \leq 2k-2$, let $D_i$ be a link diagram whose signed checkerboard graphs are the signed versions of $G_i$ and $G_{i+1}$. Then for $1 \leq i \leq 2k-3$, $D_i$ and $D_{i+1}$ share the signed checkerboard graph $G_{i+1}$. Also, $D$ and $D_1$ share the signed checkerboard graph $G=G_1$, and $D_{2k-2}$ and $D'$ share the signed checkerboard graph $G_{2k-1}=G'$. Therefore 
\[
D \sim_1 D_1 \sim_1 D_2 \sim_1 \cdots \sim_1 D_{2k-2} \sim_1 D' \text{,}
\]
and hence $D \sim_1 D'$.
\end{proof} 

We should observe that Theorem \ref{main3} remains valid if the restriction to connected checkerboard graphs is removed; however, the resulting equivalence relation is coarser than $ \sim_1 = \sim_2$. (See the earlier discussion of Figures \ref{figfour} and \ref{figfive}.) 

We should also observe that the proof of Theorem \ref{main3} makes it clear that the equivalence relation $ \sim_1 = \sim_2$ can be described in two other ways. The definitions of $\sim_1$ and $\sim_2$ in Theorem \ref{main3} could require that a checkerboard graph of $D$ and a checkerboard graph of $D'$ be related by a sign-reversing abstract duality map, or by a sign-reversing geometric duality map.

Here is an immediate consequence of Theorem \ref{main3}.

\begin{corollary}
\label{simcor}
Let $\sim'_1, \sim'_2$ be the strictest equivalence relations on links with $L \sim'_i L'$ whenever $L$ and $L'$ are links with diagrams $D$ and $D'$ such that $D \sim_i D'$. Then $\sim'_1$ and $\sim'_2$ are the same equivalence relation.
\end{corollary}

Theorem \ref{main3} and Corollary \ref{simcor} were developed as we compared several papers, which relate the signed checkerboard graphs of a link diagram to an important equivalence relation on links, called ``mutation'' \cite{CK, CG, Gre, Lip, T3}. It is not our purpose to provide detailed summaries of these papers, but we should mention that there are some potentially confusing technical differences among them. For instance, Champanerkar and Kofman \cite{CK} use the term ``mutation'' in a more restrictive sense than the other authors, as they do not allow for the Reidemeister moves; Greene \cite{Gre} only considers link diagrams in which all crossings have the same sign; and Traldi \cite{T3} allows for disconnected checkerboard graphs.

In addition to such technical differences, there is a more significant-seeming discrepancy among these references: three of them \cite{CK, CG, Gre} relate mutation to 2-isomorphism of signed checkerboard graphs, and the other two \cite{Lip, T3} relate mutation to equality of Laplacian matrices of signed checkerboard graphs. (The first knot theorist to study this Laplacian matrix was Goeritz \cite{G}, so the matrix is called the ``Goeritz matrix'' in the knot-theoretic literature.) Equality of Laplacian matrices requires isomorphism of connected, signed checkerboard graphs, because cancellation of parallel edges of opposite sign can be accomplished using Reidemeister moves. Therefore the significant-seeming discrepancy among \cite{CK, CG, Gre, Lip, T3} is essentially that \cite{CK, CG, Gre} relate mutation to 2-isomorphism of signed checkerboard graphs, while \cite{Lip, T3} relate mutation to isomorphism of signed checkerboard graphs. The ideas in the present note were developed as we tried to understand whether or not this discrepancy indicates a substantial disagreement in the knot-theoretic literature. Theorem \ref{main3} and Corollary \ref{simcor} tell us that in fact, there is not a substantial disagreement, because isomorphism and 2-isomorphism of signed checkerboard graphs define the same equivalence relation on link diagrams, and on links.

\bibliographystyle{plain}
\bibliography{dualnote}

\end{document}